# The convergence of operational Tau method for solving a class of nonlinear Fredholm fractional integro-differential equations on Legendre basis

A. Yousefi, E. Babolian, S. Javadi

Department of Computer Science, Faculty of Mathematical Sciences and Computer, Kharazmi University, Tehran, Iran

babolian@khu.ac.ir, asadyosefi@gmail.com, javadi@khu.ac.ir

**Abstract**

In this paper, we investigate approximate solutions for nonlinear Fredholm integro-differential equations (NFIDEs) of fractional order. We present an operational Tau method (OTM) by obtaining the Tau matrix representation. We solve a special class of NFIDEs based on Legendre-Tau method. The main idea in this work is to convert the fractional differential and integral parts of the desired NFIDEs to some operational matrices. By using the Sobolev inequality and some of Banach algebra properties, we prove that our proposed method converges to the exact solution in $L^2$-norm. Also, some numerical examples are presented to show the efficiency and accuracy of the proposed method.

## 1. Introduction

It is known that many phenomena in several branches of science can be described very successfully by models using mathematical tools from fractional calculus, e.g. the theory of derivatives and integrals of fractional (non-integer) order. The first definition of the fractional derivative was introduced at the end of 19-th century by Liouville and Riemann, but the concept of non-integer derivative and integral, as a generalization of the standard differential and integral calculus, was mentioned already in 1695 by Leibniz and L'Hospital. However, only in the late 1960s engineers started to be interested in this idea by the fact that descriptions of some systems are more accurate in "fractional language". Since that time fractional calculus is increasingly used to model behaviors of real systems in various fields of science and engineering. This mathematical phenomenon allows describing a real object more accurately than the classical methods, for more details see ([1], [4], [8]). In fact, during the 300 years old ago, the fractional calculus has been developed progressively up. Now, fractional calculus is applied to economics [3], dynamics of interfaces between nanoparticles and substrates [6] and image signal process [7]. Recently, several numerical methods have been proposed for solving fractional integro-differential equations (FIDEs) ([9]-[18]). Spectral methods are extremely well suited for solving FIDEs. There exist already a number of good algorithms based on the spectral Tau methods which are entirely satisfactory in terms of either efficiency or accuracy (see, for instance, [2], [11-12] and [22-25]).

Our main aim in this article is to study the equation

$$D^\alpha y(x) - \lambda \int_0^1 k(x,t)(y(t))^q dt = f(x), \quad q > 1, \tag{1}$$

subject to the initial values

$$y^{(i)}(0) = d_i, \quad i = 0,1,\ldots,m-1, \quad 1 \leq m = \lceil \alpha \rceil, \quad m \in \mathbb{N}, \qquad (2)$$

where $f \in L^2([0,1])$ and $k \in L^2([0,1] \times [0,1])$ are given functions, $y$ is the unknown function to be determined, $D^\alpha$ is the fractional derivative in the Caputo sense and q is a positive integer. Local and global uniqueness theorem of solutions NFIDEs have been shown in [19]. Vanani et al. [11] solved some kinds of the Volterra case by using operational Tau method (OTM) and presented an algorithm to find an approximate solution. Mokhtari et al. [22] apply OTM based on Legendre basis and show that is convergent for solving Volterra NFIDEs.

In the Fredholm case, Zhu et al. [18] used the second kind Chebyshev wavelet and presented an operational method for solving Eq. (1). We notice that fractional operational matrices in the case of Fredholm are different from the Volterra case. Therefore, in this article, we investigate an approximate solution by developing OTM based on shifted Legendre polynomials (for more details see [2], [5] and [22-25]) according to operational matrices in Fredholm case. We use Legendre polynomials and consider the special approximation methods based on these orthogonal polynomials, namely, the Tau approximation. Our objective is to derive simplified version of the classical algorithms. A theorem about this point is proved. We find an approximate solution which is convergent in $L^2$ −norm.

The paper is organized as follows: In section 2, some necessary definitions and mathematical tools of the fractional calculus which are required for our subsequent developments are introduced. In section 3, the OTM of the fractional derivative is obtained and proved. Section 4 is devoted to apply the OTM for solving Eq. (1) by using the shifted Legendre polynomials. After this section, we discuss about convergence analysis and then, some numerical experiments are presented in Section 6. Also, some conclusions are given in section 7.

## 2. Notation and preliminaries

Let $m \in \mathbb{N}$ be the smallest integer which is greater than or equal to $\alpha$, i.e. $m = \lceil \alpha \rceil$, the Caputo's fractional derivative operator of order $\alpha > 0$, is defined as:

$$D^\alpha y(x) = \begin{cases} J^{m-\alpha} D^m y(x), & m - 1 < \alpha \leq m, \\ D^m y(x), & \alpha = m, \end{cases} \qquad (3)$$

where

$$J^\nu y(x) = \frac{1}{\Gamma(\nu)} \int_0^x (x-t)^{\nu-1} y(t) dt, \quad \nu > 0, \quad x > 0.$$

For the Caputo's derivative we have [8]:

$$D^\alpha x^\beta = \begin{cases} 0, & \beta \in \{0,1,2,\ldots\} \text{ and } \beta < m, \\ \dfrac{\Gamma(\beta+1)}{\Gamma(\beta-\alpha+1)} x^{\beta-\alpha} & \beta \in \{0,1,2,\ldots\} \text{ and } \beta \geq m. \end{cases} \qquad (4)$$

Recall that for $\alpha \in \mathbb{N}$, the Caputo differential operator coincides with the usual differential operator. Similar to standard differentiation, Caputo's fractional differentiation is a linear operation, i.e.,

$$D^\alpha(\lambda\, g(x) + \mu\, h(x)) = \lambda D^\alpha g(x) + \mu D^\alpha h(x),$$

where $\lambda$ and $\mu$ are constants.

The Legendre polynomials $\{L_i; i = 0,1,\ldots\}$ are defined on the interval $[-1,1]$. In order to use these polynomials on the interval $[0,1]$, we define the so-called shifted Legendre polynomials by introducing the change of variable $t = 2x - 1$. Let the shifted Legendre polynomials $L_i(2x - 1)$ be denoted by $L_{1,i}(x)$, satisfying the orthogonality relation

$$\int_0^1 L_{1,i}(x)L_{1,j}(x)dx = \frac{1}{2i+1}\delta_{ij}. \tag{5}$$

The analytic form of the shifted Legendre polynomial of degree $n$ is given by

$$L_{1,i}(x) = \sum_{k=0}^{i}(-1)^{i+k}\frac{(i+k)!\,x^k}{(i-k)!\,(k!)^2}.$$

## 3. Operational Tau Method

The main objective of this section is to state the structure of OTM and its applications for the fractional calculus.

The main idea of OTM is to approximate $y \in L^2[0,1]$ with a polynomial in which $L^2[0,1]$ is the space of all functions $f: [0,1] \to \mathbb{R}$, with

$$\|f\|^2 = <f,f> := \int_0^1 f^2(x)dx < \infty.$$

Let $\Phi_x = \{\phi_i(x)\}_{i=0}^\infty$ be a sequence of arbitrary orthogonal polynomials basis. If $X_x = [1, x, x^2, \ldots]^T$, then $\Phi_x = \Phi\,X_x$. Therefore it can be easily proved that $\Phi$ is a lower triangular matrix.

The main approximate formula of the fractional derivative with any power of $y$ for solving Eq. (1) is based on the following lemma.

**Lemma 3.1 [5]** Suppose that $y(x) = \sum_{i=0}^\infty y_i\, x^i = yX_x$, then we have

$$\frac{d^r y(x)}{dx^r} = y\,M^r\,X_x, \quad r = 0,1,2,\ldots, \tag{6}$$

$$x^s y(x) = y\,N^s\,X_x, \quad s = 0,1,2,\ldots, \tag{7}$$

$$\int_a^x y(t)dt = y\,P\,X_x - y\,P\,X_a, \tag{8}$$

where $y = [y_0, y_1, y_2, \ldots]$, $X_a = [1, a, a^2, \ldots]^T$, $a \in \mathbb{R}$ and $M, N$ and $P$ are infinite matrices with only nonzero elements $M_{i+1,i} = i + 1$, $N_{i,i+1} = 1$, $P_{i,i+1} = \frac{1}{i+1}$, $i = 0,1,2,\ldots$.

**Theorem 3.2** Let $y(x) = \sum_{i=0}^\infty y_i L_{1,i}(x) = y\,\Phi\,X_x$, where $y = [y_0, y_1, y_2, \ldots]$ is the vector of unknown coefficients and $\Phi X_x$ be an orthogonal basis of polynomials in $\mathbb{R}$, then

$$D^\alpha y(x) = y\,H\,\Phi\,X_x,$$

where $H$ will be defined below.

**Proof.** By using Eq. (3) and Eq. (6), we get

$$D^\alpha y(x) = J^{m-\alpha}D^m(y\,\Phi\,X_x) = J^{m-\alpha}(y\,\Phi\,M^m\,X_x) = y\,\Phi\,M^m\,J^{m-\alpha}(X_x). \tag{9}$$

According to the definition of $X_x$ and Eq. (4), we have

$$J^{m-\alpha}(\boldsymbol{X}_x) = [J^{m-\alpha}(1), J^{m-\alpha}(x), J^{m-\alpha}(x^2), \dots]^T$$
$$= \left[\frac{\Gamma(1)x^{m-\alpha}}{\Gamma(m-\alpha+1)}, \frac{\Gamma(2)x^{m-\alpha+1}}{\Gamma(m-\alpha+2)}, \frac{\Gamma(3)x^{m-\alpha+2}}{\Gamma(m-\alpha+3)}, \dots\right]^T = \boldsymbol{\Gamma}\,\Lambda, \tag{10}$$

where $\boldsymbol{\Gamma}$ is an infinite diagonal matrix with elements

$$\Gamma_{i,i} = \frac{\Gamma(i+1)}{\Gamma(m-\alpha+i+1)}, \qquad i = 0,1,\dots,$$

and

$$\Lambda = [x^{m-\alpha}, x^{m-\alpha+1}, x^{m-\alpha+2}, \dots]^T.$$

Let us approximate $x^{m-\alpha+j}$ as

$$x^{m-\alpha+j} = \sum_{i=0}^{\infty} a_{ji} L_{1,i}(x) = \boldsymbol{a}_j\,\boldsymbol{\Phi}\,\boldsymbol{X}_x, \quad \boldsymbol{a}_j = [a_{j0}, a_{j1}, a_{j2}, \dots].$$

Therefore, we obtain

$$\Lambda = [\boldsymbol{a}_0\,\boldsymbol{\Phi}\,\boldsymbol{X}_x, \boldsymbol{a}_1\,\boldsymbol{\Phi}\,\boldsymbol{X}_x, \boldsymbol{a}_2\,\boldsymbol{\Phi}\,\boldsymbol{X}_x, \dots]^T = \boldsymbol{A}\,\boldsymbol{\Phi}\,\boldsymbol{X}_x, \tag{11}$$

$$\boldsymbol{A} = [\boldsymbol{a}_0, \boldsymbol{a}_1, \boldsymbol{a}_2, \dots]^T.$$

Substituting Eq. (11) in Eq. (10) and then Eq. (9), we obtain:

$$D^\alpha y(x) = \boldsymbol{y}\,\boldsymbol{\Phi}\,\boldsymbol{M}^m\,\boldsymbol{\Gamma}\,\Lambda = \boldsymbol{y}\,\boldsymbol{\Phi}\,\boldsymbol{M}^m\,\boldsymbol{\Gamma}\,\boldsymbol{A}\,\boldsymbol{\Phi}\,\boldsymbol{X}_x = \boldsymbol{y}\,\boldsymbol{H}\,\boldsymbol{\Phi}\,\boldsymbol{X}_x, \quad \boldsymbol{H} = \boldsymbol{\Phi}\,\boldsymbol{M}^m\,\boldsymbol{\Gamma}\,\boldsymbol{A}. \tag{12}$$

## 4. Applications and Algorithm of Our Approach

Consider the nonlinear Fredholm fractional integro-differential equation of the second kind given by Eq. (1). In this section, we introduce the operational Tau formula for Eq. (1) by using the shifted Legendre polynomials.

For this purpose, we state the following lemma.

**Lemma 4.1 [11].** Let $\boldsymbol{X}_x, \boldsymbol{y}$ and $\boldsymbol{\Phi}$ be the vectors and matrix defined in the previous section. Then, we have:

$$\boldsymbol{X}_x\,\boldsymbol{y}\,\boldsymbol{\Phi}\,\boldsymbol{X}_x = \boldsymbol{Y}\,\boldsymbol{X}_x, \tag{13}$$

where $\boldsymbol{Y}$ is an upper triangular matrix with entries

$$Y_{i,j} = \begin{cases} \sum_{k=0}^{\infty} y_k \phi_{k,j-1}, & j \geq i, \ i,j = 0,1,\dots, \\ 0, & j < i, \ i,j = 0,1,\dots. \end{cases}$$

In addition, if we suppose $y(x) = \boldsymbol{y}\,\boldsymbol{\Phi}\,\boldsymbol{X}_x$ represents a polynomial, then for any positive integer $p$, the following relation is valid

$$y^q(x) = \boldsymbol{y}\,\boldsymbol{\Phi}\,\boldsymbol{Y}^{q-1}\,\boldsymbol{X}_x.$$

The above relation can be used to linearize the function $F(x,t,y(t)) = k(x,t)y^q(t)$. Hence we have

$$F(x,t,y(t)) = k(x,t)y^q(t) = (\boldsymbol{X}_x^T\,\boldsymbol{K}\,\boldsymbol{X}_t)(\boldsymbol{y}\,\boldsymbol{\Phi}\,\boldsymbol{Y}^{q-1}\,\boldsymbol{X}_x)$$
$$= \boldsymbol{X}_t^T\,\boldsymbol{K}^T\,(\boldsymbol{X}_x\,\boldsymbol{y}\,\boldsymbol{\Phi}\,\boldsymbol{Y}^{q-1}\,\boldsymbol{X}_x). \tag{14}$$

For simplicity, we define $\Lambda := \Phi . Y^{q-1}$, i.e.
$$\mathbf{B} = [\mathbf{B}_0 | \mathbf{B}_1 | \mathbf{B}_2 | \dots ] = [\Phi_0 | \Phi_1 | \Phi_2 | \dots] \times \left[ Y_0^{(q-1)} \Big| Y_1^{(q-1)} \Big| Y_2^{(q-1)} \Big| \dots \right].$$

Thus, the parenthesis on the right Eq. (14) can be written as
$$X_x \, y \, \Phi \, Y^{q-1} \, X_x = X_x \, [y \, \mathbf{B}_0 | y \, \mathbf{B}_1 | y \, \mathbf{B}_2 | \dots ] \, X_x$$
$$= \begin{bmatrix} y\mathbf{B}_0 & y\mathbf{B}_1 & \cdots \\ y\mathbf{B}_0 x & y\mathbf{B}_1 x & \cdots \\ y\mathbf{B}_0 x^2 & y\mathbf{B}_1 x^2 & \cdots \\ & \vdots & \end{bmatrix} \begin{bmatrix} 1 \\ x \\ x^2 \\ \vdots \end{bmatrix}$$
$$= \begin{bmatrix} \sum_{i=0}^{\infty} y\mathbf{B}_i \, x^i \\ \sum_{i=0}^{\infty} y\mathbf{B}_i \, x^{i+1} \\ \sum_{i=0}^{\infty} y\mathbf{B}_i \, x^{i+2} \\ \vdots \end{bmatrix} = \begin{bmatrix} y\mathbf{B}_0 & y\mathbf{B}_1 & \cdots \\ 0 & y\mathbf{B}_0 & \cdots \\ 0 & 0 & \cdots \\ & \vdots & \end{bmatrix} \begin{bmatrix} 1 \\ x \\ x^2 \\ \vdots \end{bmatrix}$$
$$= \Delta X_x, \tag{15}$$

where $\Delta$ is an upper triangular matrix with entries
$$\Delta_{i,j} = \begin{cases} y\mathbf{B}_{j-i}, & j \geq i, \quad i, j = 0, 1, \dots, \\ 0, & j < i, \quad i, j = 0, 1, \dots . \end{cases}$$

By substituting Eq. (15) to Eq. (14), we have
$$F(x, t, y(t)) = X_t^T \, K^T \, \Delta X_x. \tag{16}$$

By virtue of Eq. (16) and Eq. (12), we approximate Eq. (1) by the following operational form:

$$y \, H \, \Phi \, X_x - \lambda P^T \, K^T \, \Delta X_x = f \, X_x, \tag{17}$$

where $P^T = \int_0^1 X_t^T \, dt$.

So, the residual of Eq. (1) can be written as
$$R(x) = \{y H \Phi - \lambda P^T K^T \Delta - f\} X_x =: RX_x. \tag{18}$$

The exact solution of Eq. (1) is obtained when $R(x) \equiv 0$. Since the shifted Legendre polynomials are complete, this is achievable by solving the following infinite system of algebraic equations
$$(R, L_{1,k}) = 0, \quad k = 0, 1, \dots . \tag{19}$$

In addition, by using Eq. (6), the initial conditions can be written as
$$y M^i X_x \big|_{x=0} = d_i, \quad i = 0, 1, \dots, m - 1. \tag{20}$$

In practice, we truncate the series $y(x) = \sum_{i=0}^{\infty} y_i L_{1,i}(x) = y \, \Phi \, X_x$ to finite number of terms. Hence choosing $N + 1 - m$ of the first equations in Eq. (19) along with Eq. (20) lead to the following system
$$(R, L_{1,k}) = 0, \quad k = 0, 1, \dots, N - m,$$
$$y \, M^i \, X_0 = d_i, \quad i = 0, 1, \dots, m - 1.$$

Solving the above system yields the unknown vector $\boldsymbol{y} = [y_0, y_1, \ldots, y_N]$.

## 5. Convergence Analysis of the Legendre Tau Method

In this section, we present a general approach to the convergence analysis for nonlinear Fredholm fractional integro-differential in $L^2$ −norm. Here, there are some properties that can be found in [20,25].

We first give the definition of Sobolev space on $I = [0,1]$: $H^m(I)$ as the space of all functions $y$ on $I$ such that $y$ and all its weak derivatives up to order $m$ are in $L^2(I)$. $H^m(I)$ is a Hilbert space with the norm and semi norm defined respectively by

$$\|y\|^2_{H^m(I)} = \sum_{i=0}^{m} \|y^{(i)}\|^2_{L^2(I)}, \tag{21}$$

and

$$|y|^2_{H^{l:N}(I)} = \sum_{i=\min(l:N)}^{N} \|y^{(i)}\|^2_{L^2(I)}. \tag{22}$$

Let $P_N(I)$ be the space of polynomials with degree$\leq N$ on $I$ and $p_N$ be the orthogonal projection operator from $L^2(I)$ upon $P_N(I)$. Then, for any function $f$ in $L^2(I)$, $p_N f$ belongs to $P_N(I)$ and satisfies

$$\int_0^1 (f - p_N f)(t) g(t) dt = 0, \quad \forall g \in P_N(I).$$

Also, the following relations with shifted Legendre polynomials may readily be obtained for $1 \leq l, k \leq m$:

$$\|y - p_N(y)\|_{H^l(I)} \leq C_1 N^{2l - \frac{1}{2} - k} |y|_{H^{k:N}(I)}, \tag{23}$$

$$\|y - p_N(y)\|_{L^2(I)} \leq C_2 N^{-k} |y|_{H^{k:N}(I)}, \tag{24}$$

where $y \in H^m(I)$, and $C_1$ and $C_2$ are constants independent of $N$.

It is trivial that, for any Banach algebra space $V$, we can write [21]:

$$\forall x, y \in V, \quad \|xy\| \leq \|x\| \|y\|,$$
$$\|x^n\| \leq \|x\|^n.$$

At last, we apply the Sobolev inequality.

**Lemma 5.1 [20]** (Sobolev inequality) Assume $(a, b)$ be a bounded interval of the real line. For any function $y \in H^1(a, b)$, the following inequality holds:

$$\|y\|_{L^\infty} \leq \left(\frac{1}{b-a} + 2\right)^{\frac{1}{2}} \|y\|^{\frac{1}{2}}_{L^2(a,b)} \|y\|^{\frac{1}{2}}_{H^1(a,b)}.$$

**Lemma 5.2[4]** If $m - 1 < \alpha \leq m, m \in \mathbb{N}$ and $f \in C^m(0,1)$, then

$$D^\alpha J^\alpha f(x) = f(x),$$

and

$$J^\alpha D^\alpha f(x) = f(x) - \sum_{k=0}^{m-1} f^{(k)}(0^+)\frac{x^k}{k!}, \quad x > 0. \tag{25}$$

**Lemma 5.3** The fractional integro-differential equation (1) together (2) is equivalent to the following nonlinear integral equation

$$y(x) = \sum_{k=0}^{m-1} d_k \frac{x^k}{k!} + \frac{1}{\Gamma(\alpha)} \int_0^x (x-s)^{\alpha-1} f(x) ds$$

$$+ \frac{\lambda}{\Gamma(\alpha)} \int_0^x (x-s)^{\alpha-1} \int_0^1 k(x,t) y^q(t) \, ds. \tag{26}$$

**Proof.** It can be proved by applying the fractional integral operator to both sides of Eq. (1), and using initial values and **Lemma 5.2** to get Eq. (26).

Now, we shall prove the main result of our paper in this section. In the following theorem, an error estimation for the approximate solution of Eq. (1) with supplementary conditions of Eq. (2) is obtained. Let $e_N(x) = y(x) - y_N(x)$, be the error function of the approximate solution $y_N$ obtained by the spectral Legendre Tau method for $y(x)$.

**Theorem 5.2** For sufficiently large $N$, the spectral Legendre Tau approximation converges to the exact solution in $L^2$-norm, i.e.

$$\|e_N\|_{L^2(I)} = \|y - y_N\|_{L^2(I)} \to 0, \quad as\ N \to \infty, \quad I = (0,1).$$

Proof. According to Legendre Tau method, we have

$$y_N(x) = \sum_{k=0}^{m-1} d_k \frac{x^k}{k!} + \frac{1}{\Gamma(\alpha)} p_N \left( \int_0^x (x-s)^{\alpha-1} f(x) \, ds \right)$$

$$+ \frac{\lambda}{\Gamma(\alpha)} p_{N,N} \left( \int_0^x (x-s)^{\alpha-1} \int_0^1 k(x,t) \, y_N^q(t) \, dt \, ds \right). \tag{27}$$

By subtracting Eq. (27) from Eq. (26) by this assumption, we have

$$e_N(x) = \frac{1}{\Gamma(\alpha)} \left( \int_0^x (x-s)^{\alpha-1} f(x) ds - p_N \left( \int_0^x (x-s)^{\alpha-1} f(x) \, ds \right) \right)$$

$$+ \frac{\lambda}{\Gamma(\alpha)} \left( \int_0^x (x-s)^{\alpha-1} \int_0^1 k(.,t) y^q(t) \, ds \right.$$

$$\left. - p_N \left( \int_0^x (x-s)^{\alpha-1} \int_0^1 k(.,t) \, y_N^q(t) \, dt \, ds \right) \right). \tag{28}$$

By taking the $L^2$-norm, we can get

$$\|e_N\|_{L^2(I)} \leq \left\|\frac{1}{\Gamma(\alpha)}\left(\int_0^x (x-s)^{\alpha-1}f(x)ds - p_N\left(\int_0^x (x-s)^{\alpha-1}f(x)\,ds\right)\right)\right\|_{L^2(I)}$$

$$+ \left\|\frac{\lambda}{\Gamma(\alpha)}\left(\int_0^x (x-s)^{\alpha-1}\int_0^1 k(x,t)y^q(t)\,dt\,ds\right.\right.$$

$$\left.\left. - p_N\left(\int_0^x (x-s)^{\alpha-1}\int_0^1 k(x,t)\,y_N^q(t)\,dt\,ds\right)\right)\right\|_{L^2(I)}$$

$$= \left\|\frac{1}{\Gamma(\alpha)}\left(\int_0^x (x-s)^{\alpha-1}f(x)ds - p_N\left(\int_0^x (x-s)^{\alpha-1}f(x)\,ds\right)\right)\right\|_{L^2(I)}$$

$$+ \left\|\frac{\lambda}{\Gamma(\alpha)}\left(\int_0^x (x-s)^{\alpha-1}\int_0^1 k(.,t)y^q(t)\,dt\,ds\right.\right.$$

$$- \int_0^x (x-s)^{\alpha-1}\int_0^1 k(.,t)y_N^q(t)\,dt\,ds$$

$$+ \int_0^x (x-s)^{\alpha-1}\int_0^1 k(.,t)y_N^q(t)\,dt\,ds$$

$$\left.\left. - p_N\left(\int_0^x (x-s)^{\alpha-1}\int_0^1 k(.,t)\,y_N^q(t)\,dt\,ds\right)\right)\right\|_{L^2(I)} \tag{29}$$

From (24), we may write

$$\|e_N\|_{L^2(I)} \leq C_2\, N^{-k}\left|\int_0^x (x-s)^{\alpha-1}f(x)ds\right|_{H^{k:N}(I)}$$

$$+ \left\|\int_0^x (x-s)^{\alpha-1}\int_0^1 k(.,t)\left(y^q(t)-y_N^q(t)\right)dt\,ds\right\|_{L^2(I)}$$

$$+ C_2\, N^{-k}\left|\int_0^x (x-s)^{\alpha-1}\int_0^1 k(.,t)y_N^q(t)\,dt\,ds\right|_{H^{k:N}(I)}$$

$$\leq C_3 N^{-k}\|f\|_\infty + C_4 \max_{x,t\in I}|k(x,t)|\,\|y^q - y_N^q\|_\infty + C_5\, N^{-k}\max_{x,t\in I}|k(x,t)|\,\|y_N^q\|_\infty$$

$$\leq C_3 N^{-k}\|f\|_\infty + C_4 \max_{x,t\in I}|k(x,t)|\,\|y^q - y_N^q\|_\infty$$

$$+ C_5\, N^{-k}\max_{x,t\in I}|k(x,t)|\left(\|y^q - y_N^q\|_\infty + \|y^q\|_\infty\right)$$

$$= E_1 + E_2 + E_3 \tag{30}$$

where $C_3, C_4$ and $C_4$ are constants independent of $N$ and

$$E_1 = C_3 N^{-k}\|f\|_\infty,$$
$$E_2 = C_4 \max_{x,t\in I}|k(x,t)|\,\|y^q - y_N^q\|_\infty,$$

$$E_3 = C_5 \, N^{-k} \max_{x,t \in I} |k(x,t)| \left( \|y^q - y_N^q\|_\infty + \|y^q\|_\infty \right). \tag{31}$$

Since $k \geq 1$ and
$$N^{-k} \to 0 \text{ as } N \to \infty \text{ for all } k \geq 1,$$
we can conclude $E_1$ and $E_2$ tend to zero for sufficiently large $N$. The rest of our proof is to show that
$$\|y^q - y_N^q\|_\infty \to 0 \quad \text{as } N \to \infty.$$

To do this, we may use binomial formula as follow
$$y^q(x) - y_N^q(x) = (y(x) - y_N(x))\left(\sum_{i=0}^{q-1} y^{q-i} \, y_N^i\right)$$
$$= e_N \left( \sum_{i=0}^{q-1} y^{q-i} (y - e_N)^i \right)$$
$$= e_N \left( \sum_{i=0}^{q-1} y^{q-i} \sum_{j=0}^{i} \binom{i}{j} y^j (-e_N)^{i-j} \right)$$
$$= e_N \left( \sum_{i=0}^{q-1} \sum_{j=0}^{i} \binom{i}{j} y^{q-i+j} (-e_N)^{i-j} \right). \tag{32}$$

By virtue of Banach algebra theory, we have
$$\|y^q - y_N^q\|_\infty \leq \|e_N\| \times \sum_{i=0}^{q-1} \sum_{j=0}^{i} \binom{i}{j} \|e_N\|_\infty^{i-j} \|y\|_\infty^{q+j-i}. \tag{33}$$

Finally, it is sufficient to prove that $\|e_N\|_\infty \to 0$ as $N \to \infty$.
On the other hand, from Lemma 5.1, i.e. Sobolev inequality together with (23) and (24), we conclude that $\|e_N\|_\infty \to 0$ as $N \to \infty$, because
$$\|e_N\|_\infty \leq 3^{\frac{1}{2}} \|e_N\|_{L^2(I)}^{\frac{1}{2}} \|e_N\|_{H^1(I)}^{\frac{1}{2}}$$
$$\leq 3^{\frac{1}{2}} \left( C_2 \, N^{-k} |y|_{H^{k:N}(I)} \right)^{\frac{1}{2}} \left( C_1 \, N^{\frac{3}{2}-k} |y|_{H^{k:N}(I)} \right)^{\frac{1}{2}}$$
$$= C_6 \, N^{\frac{3}{4}-k} |y|_{H^{k:N}(I)}, \tag{34}$$
where
$$C_6 = \sqrt{3} \, C_1 C_2.$$

Since $k \geq 1$, then we can easily conclude $N^{\frac{3}{4}-k} \to 0$ as $N \to \infty$ and therefore,
$$\|e_N\|_\infty \to 0 \text{ as } N \to \infty.$$

Similarly, $\|e_N\| \times \|e_N\|_\infty^{i-j}$ for $i \geq j$, tends to zero as $N \to \infty$. Then
$$E_2 \to 0 \text{ as } N \to \infty.$$

Hence the theorem is proved.

## 6. Numerical example

To show efficiency of our numerical method, the following examples are considered.

**Example 1.** Consider the second kind nonlinear Fredholm fractional integro-differential equation [18]

$$D^{\frac{1}{2}}y(x) - \int_0^1 xt(y(t))^4 \, dt = \frac{1}{\Gamma\left(\frac{1}{2}\right)}\left(\frac{8}{3}x^{\frac{3}{2}} - 2x^{\frac{1}{2}}\right) - \frac{x}{1260}, \quad 0 \le x < 1,$$

with the initial condition: $y(0) = 0$, and exact solution $y(x) = x^2 - x$.
We have solved this example using OTM with shifted Legendre polynomials and approximations are obtained as follows:

$$\begin{aligned} n = 0: & \quad y_0(x) = 0, \\ n = 1: & \quad y_1(x) = -0.97238662853x, \\ n = 2: & \quad y_2(x) = x^2 - x, \\ n = 3: & \quad y_3(x) = x^2 - x, \end{aligned}$$

and so on. Therefore, we obtain $y(x) = x^2 - x$ which is the exact solution of the problem.

**Example 2.** Consider the following NFIDEs [18]

$$D^{\frac{5}{3}}y(x) - \int_0^1 (x+t)^2 (y(t))^3 \, dt = \frac{6}{\Gamma\left(\frac{1}{3}\right)} x^{\frac{1}{3}} - \frac{x^2}{7} - \frac{x}{4} - \frac{1}{9}, \quad 0 \le x < 1,$$

with the initial conditions: $y(0) = y'(0) = 0$.
The exact solution is $y(x) = x^2$, and similarly, we obtain the exact solution of this problem for n = 2,3, ....

**Example 3.** Consider the following second kind NFIDEs [18]

$$D^{\alpha}y(x) - \int_0^1 xt(y(t))^2 \, dt = 1 - \frac{x}{4}, \quad 0 \le x < 1, 0 \le \alpha < 1,$$

subject to the initial condition $y(0) = 0$. The only case which we know the exact solution is $\alpha = 1$ with $y(x) = x$.
This problem have been solved for $n = 10$ with different $\alpha$ by OTM and the results have been compared with the solution obtained by the second kind Chebyshev wavelet [18].
For the case $\alpha = 1$, we obtain $y_1(x) = x$, that is the exact solution.
In Table 1, the results show the accuracy of the proposed method. Furthermore, the OTM in **Figure 1** is compared with the method in [18] for some $0 < \alpha \le 1$.
In the theory of fractional calculus [8], it is clear that when the fractional derivative $\alpha$ $(m-1 < \alpha \le m)$ tends to $m$ $(m \in \mathbb{N})$, then the approximating solution continuously tends to the exact solution of the problem with $m = \lceil \alpha \rceil$. A closer look at Table 1, confirms this fact. Also, we compare the graphs of exact solutions and the approximating solutions obtained by the Chebyshev wavelet method in figure 1.

Table 1. The approximation solution of Example 3 for some $0 < \alpha \le 1$.

| T | OTM for $\alpha=\frac{1}{4}$ | OTM for $\alpha=\frac{1}{2}$ | OTM for $\alpha=\frac{3}{4}$ | OTM for $\alpha=1$ |
|---|---|---|---|---|
| ۰ | 0.0000000000 | 0.0000000000 | 0.0000000000 | 0.0000000000 |
| ۰,۱ | 0.2369652099 | 0.1650204533 | 0.1247524923 | 0.1000000000 |
| ٠,٢ | 0.4528141232 | 0.3197389836 | 0.2455477221 | 0.2000000000 |
| ٠,٣ | 0.6475467400 | 0.4641555910 | 0.3623856895 | 0.3000000000 |
| ٠,٤ | 0.8211630604 | 0.5982702754 | 0.4752663945 | 0.4000000000 |
| ٠,٥ | 0.9736630842 | 0.7220830368 | 0.5841898369 | 0.5000000000 |
| ٠,٦ | 1.1050468117 | 0.8355938753 | 0.6891560170 | 0.6000000000 |
| ٠,٧ | 1.2153142426 | 0.9388027908 | 0.7901649346 | 0.7000000000 |
| ٠,٨ | 1.3044653770 | 1.0317097834 | 0.8872165897 | 0.8000000000 |
| ٠,٩ | 1.3725002150 | 1.1143148529 | 0.9803109824 | 0.9000000000 |
| ١ | 1.4194187565 | 1.1866179995 | 1.0694481126 | 1.0000000000 |

**Figure 1.** Plots of approximate solution obtained by our method and the Chebyshev wavelet method for Example 3 for $\alpha = \frac{1}{4}, \frac{1}{2}, \frac{3}{4}, \alpha = 1$, respectively.

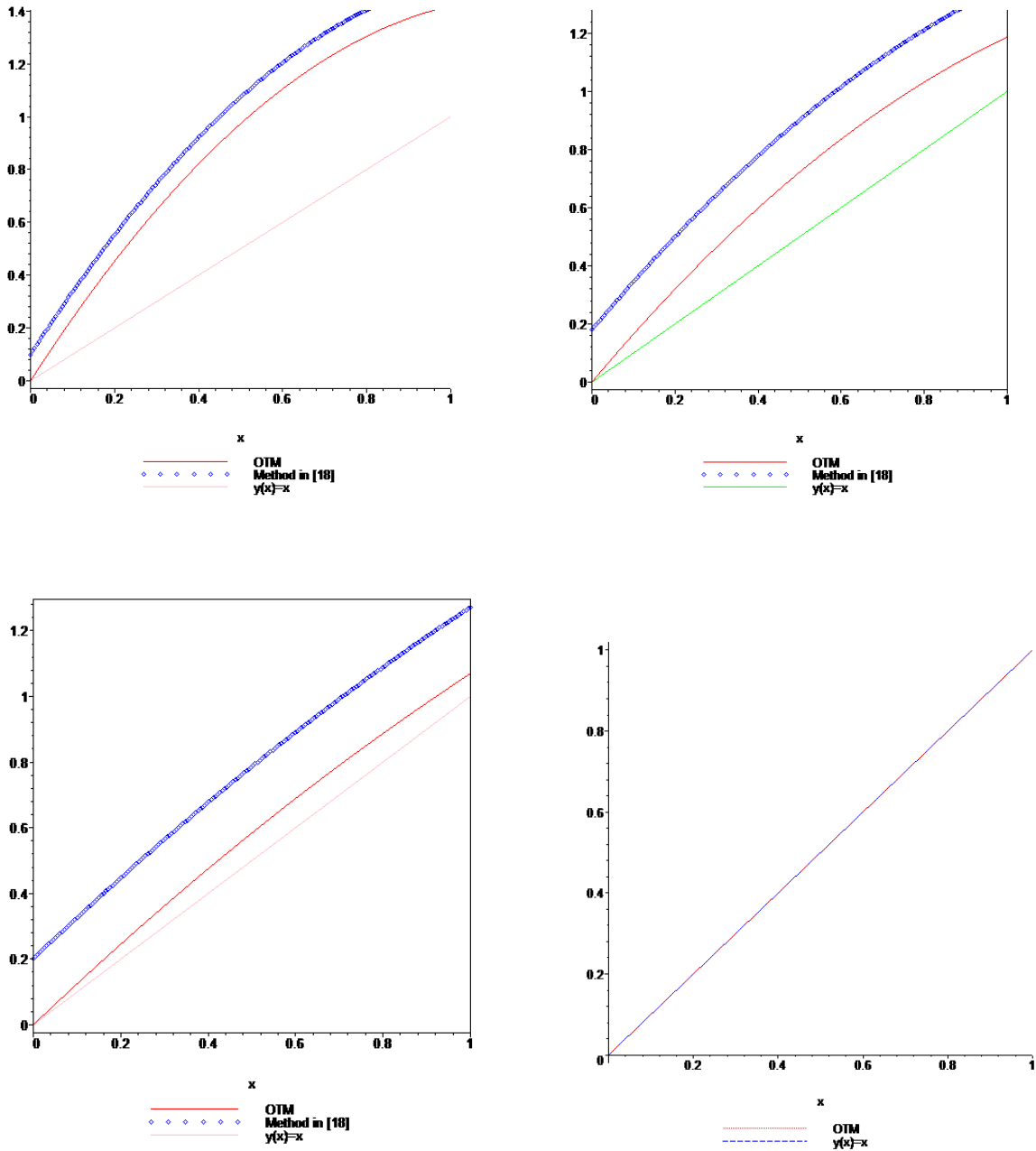

The approximate solution by OTM is shown by line, and the approximate solution obtained by Zhu[18] by point. This problem has been integrated over the interval $I = [0,1]$ such that maximal errors have been evaluated with various number $N$ of steps. **Table 2** and **Figure 2** show the exponential rate of convergence predicted by our method, i.e. Legendre-Tau method.

Table 2. The OTM errors of Example 3 for some $\alpha = \frac{1}{2}$

| N | Maximal error in Legendre-Tau method |
|---|---|
| 2 | $1.23 \times 10^{-2}$ |
| 4 | $1.56 \times 10^{-4}$ |
| 8 | $3.54 \times 10^{-9}$ |
| 16 | $1.77 \times 10^{-15}$ |
| 32 | $2.09 \times 10^{-33}$ |

**Figure 2.** The Legendre Tau approximation errors of example 3 for $N = 8$ and $16$ and $\alpha = \frac{1}{2}$

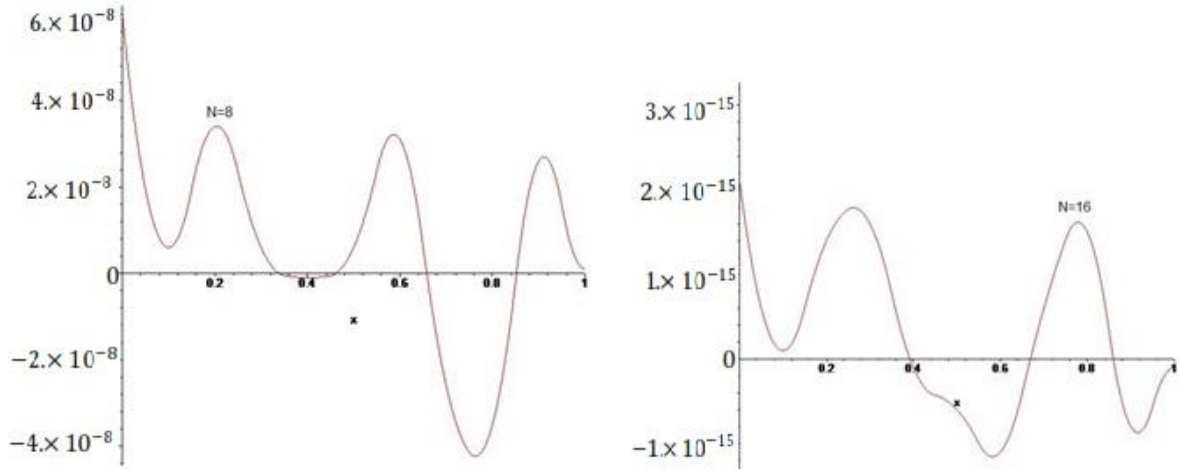

**Example 4.** Consider the following NFIDEs

$$D^{\frac{3}{2}}y(x) - \int_0^1 e^{xt}\bigl(y(t)\bigr)^2 \, dt = erf(x)e^x - \frac{e^{x+3} - 1}{x + 3}, \quad 0 \le x \le 1,$$

subject to the initial condition $y(0) = 1$. $erf$ is the error function defined for all complex $x$ by the following relation

$$\operatorname{erf}(x) = \frac{2}{\pi} \int_0^x e^{-t^2} \, dt.$$

$y(x) = e^x$ is the exact solution. **Table 3** shows the exact and approximate solution for $N = 16, 32$ and $64$. **Figure 3** illustrate the maximal absolute error with $N = 16, 32$ and $32$. From **Table 3** and **Figure 3** we can see the approximate solutions are in a good agreement with the exact solutions.

Table 3. Exact and approximation solution of Example 4

| T | Exact solution | OTM for N=8 | OTM for N=32 | OTM for $N=64$ |
|---|---|---|---|---|
| 0 | 1.000000000 | 1.000000000 | 1.000000000 | 1.000000000 |
| 0.1 | 1.105170918 | 1.105170903 | 1.105170918 | 1.105170918 |
| 0.2 | 1.221402758 | 1.221402742 | 1.221402758 | 1.221402758 |
| 0.3 | 1.349858808 | 1.349858812 | 1.349858808 | 1.349858808 |
| 0.4 | 1.491824698 | 1.491824665 | 1.491824698 | 1.491824698 |
| 0.5 | 1.648721271 | 1.648721224 | 1.648721271 | 1.648721271 |
| 0.6 | 1.822118800 | 1.822118824 | 1.822118800 | 1.822118800 |
| 0.7 | 2.013752707 | 2.013752713 | 2.013752707 | 2.013752707 |
| 0.8 | 2.225540928 | 2.225540925 | 2.225540928 | 2.225540928 |
| 0.9 | 2.459603111 | 2.459603134 | 2.459603111 | 2.459603111 |
| 1 | 2.718281828 | 2.718281859 | 2.718281828 | 2.718281828 |

Figure 3. The Legendre Tau approximation errors of example 4 for $N=8$

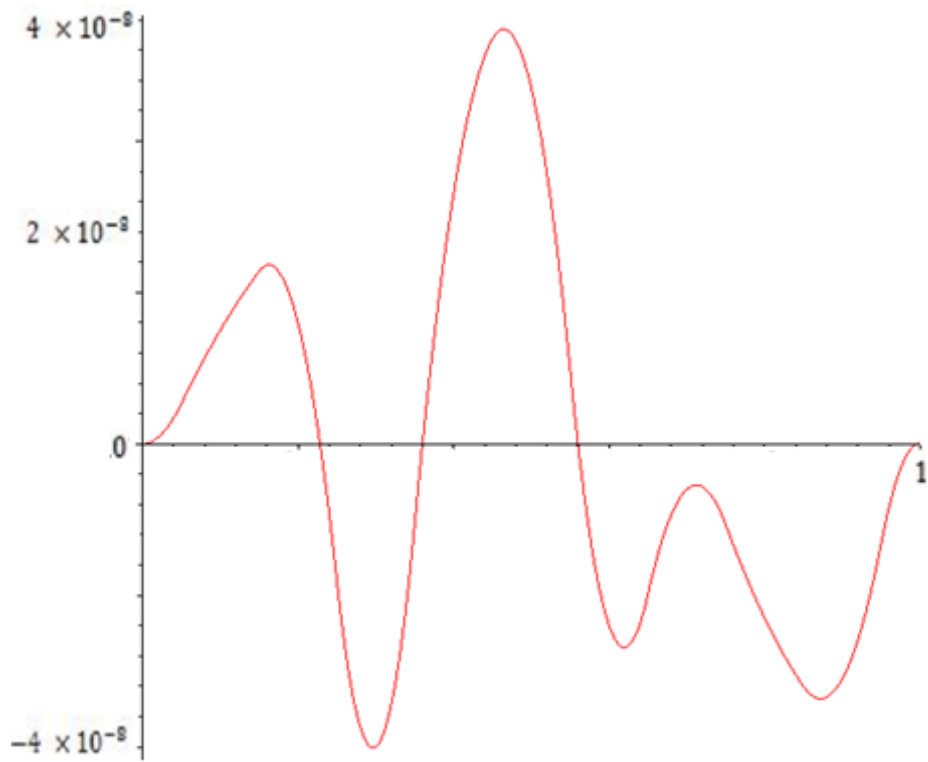

## 7. Conclusion

We applied a general formulation of OTM for fractional derivatives, which is used to approximate the numerical solution of a class of NFIDEs. Our approach was based on the shifted Legendre polynomials. The main advantage of the present method (OTM) is its simplicity, and is more convenient for computer algorithms. Furthermore, this method of solution yields the desired accuracy under certain conditions.

Orthogonal bases such as Legendre polynomials are used to reduce the run time of the source programs and the computation time due to their orthogonal properties. All these advantages of the OTM to solve nonlinear problems assert the method as a convenient, reliable, effective and powerful tool.

The fractional derivatives are described in the Caputo sense, because the Caputo fractional derivative allows traditional initial conditions to be included in the formulation of the problem. Our results given in the previous section demonstrate the good accuracy of Legendre-Tau method in fractional problems in the Caputo sense. Moreover, in the special case, only small number of shifted Legendre polynomials is needed to obtain a satisfactory result.